\documentclass[journal]{IEEEtran}

\usepackage{amsmath,amssymb,bm,amsthm}
\usepackage{amsfonts}
\usepackage{mathtools}
\usepackage{subfig}
\usepackage{graphicx}
\usepackage{booktabs}
\usepackage[english]{babel}
\usepackage{float}
\usepackage{cite}
\DeclarePairedDelimiter{\ceil}{\lceil}{\rceil}

\usepackage{comment,physics}
\usepackage[margin=1in]{geometry}
\usepackage{stackengine}
\usepackage{lineno}
\usepackage{xcolor}
\usepackage{array,multirow}
\usepackage{makecell}

\usepackage{scalerel}
\usepackage{todonotes}

\def\dashint{\,\ThisStyle{\ensurestackMath{%
\stackinset{c}{.2\LMpt}{c}{.5\LMpt}{\SavedStyle-}{\SavedStyle\phantom{\int}}}%
        \setbox0=\hbox{$\SavedStyle\int\,$}\kern-\wd0}\int}

\usepackage[linesnumbered,ruled,norelsize]{algorithm2e}

\makeatletter
\newcommand{\removelatexerror}{\let\@latex@error\@gobble}
\makeatother

\begin{document}
\bstctlcite{IEEEexample:BSTcontrol}
    \title{A Complete Helmholtz Decomposition on Multiply Connected Subdivision Surfaces and Its Application to Integral Equations}
  \author{A.~M.~A.~Alsnayyan,~\IEEEmembership{Student Member,~IEEE,} L.~Kempel,~\IEEEmembership{Fellow,~IEEE,}
      and~B.~Shanker,~\IEEEmembership{Fellow,~IEEE}

  \thanks{The  authors  acknowledge  computing  support  from  the  HPC  Center  at Michigan  State  University,  financial  support  from  NSF  via  CMMI-1725278.}
  
  \thanks{The  authors  are  with  the  Department  of  Electrical  and  Computer Engineering, Michigan State University, East Lansing, MI 48824-1226 USA. (e-mail: alsnayy1@msu.edu).}
}

\markboth{IEEE TRANSACTIONS ON ANTENNAS AND PROPAGATION
}{Alsnayyan \MakeLowercase{\textit{et al.}}: Well conditioned Integral Equations}

\maketitle

\begin{abstract}

The analysis of electromagnetic scattering in the isogeometric analysis (IGA) framework based on Loop subdivision has long been restricted to simply-connected geometries. The inability to analyze multiply-connected objects is a glaring omission. In this paper, we address this challenge. IGA provides seamless integration between the geometry and analysis by using the same basis set to represent both. In particular, IGA methods using subdivision basis sets exploit the fact that the basis functions used for surface description are smooth (with continuous second derivatives) almost everywhere. On simply-connected surfaces, this permits the definition of basis sets that are divergence-free and curl-free. What is missing from this suite is a basis set that is both divergence-free \emph{and} curl-free, a necessary ingredient for a complete Helmholtz decomposition of currents on multiply-connected structures. In this paper, we effect this missing ingredient numerically using random polynomial vector fields. We show that this basis set is analytically divergence-free and curl-free. Furthermore, we show that these basis recovers curl-free, divergence-free, and curl-free and divergence-free fields. Finally, we use this basis set to discretize a well-conditioned integral equation for analyzing perfectly conducting objects and demonstrate excellent agreement with other methods. 
\end{abstract}


\IEEEpeerreviewmaketitle


\section{Introduction}

\IEEEPARstart{S}{urface} integral equation based solvers have become a workhorse for electromagnetic (EM) analysis for problems ranging from modeling and simulation of coupled circuit–EM problems \cite{swanson2003microwave} to device optimization \cite{sideris2019ultrafast,garza2020boundary} to scattering and radiation \cite{peterson1998computational}. Given the range of applications, a number of different types of equations have been developed for analysis of composite objects \cite{peterson1998computational,vico2016decoupled,liu2018potential}, and within the past two decades there has been an extensive body of work that has examined fundamental nuances of integral equations. These include developing methods to overcome critical bottlenecks, including understanding the nuances and ramifications of discretizing these equations. This extends to various efforts to understand low frequency breakdown \cite{hsiao1997mathematical,calderon_EFIE}, develop well conditioned formulations \cite{adrian2021electromagnetic,boubendir2014well,vico2016decoupled,liu2018potential}, investigate accuracy and convergence \cite{warnick2008numerical}, and so on. We note that most of this analysis has been applied to geometric models that are Lagrangian, i.e., tessellation with discontinuous normals.  

A more recent trend has been the development of  isogeometric analysis (IGA) methods, or methods wherein the same basis set is used to represent both the geometry and the physics on the geometry \cite{takahashi2022isogeometric,wolf2021isogeometric,Jie,simpson2018isogeometric,buffa2010isogeometric}. The motivation is to permit analysis directly on computer aided design (CAD) models as opposed to a meshed model. Typically, CAD models are higher order smooth throughout the geometry (other than where necessary) and have significantly fewer control nodes. The ability to define basis sets directly on the CAD model offers a number of advantages. Typical CAD models fall under two categories, those either using non-uniform B-splines (NURBS) or subdivision. While IGA methods for electromagnetic analysis on NURBS exists \cite{buffa2010isogeometric,takahashi2022isogeometric,dolz2019isogeometric}, the challenge with using NURBS to model geometries is that they may not be watertight and are sometimes discontinuous. These hinder application to  complex geometries \cite{Bazilevs2010}, ans was one of the key motivating factors for developing subdivision surfaces \cite{zorin2007subdivision}.

Subdivision surfaces (especially Loop subdivision) has been extremely popular in the computer graphics industry due to the ease with which one can represent complex topologies, its scalability, inherently multiresolution features, efficiency and ease of implementation; it overcomes several drawbacks of NURBS models \cite{zorin2007subdivision}. More importantly, the surface representation is $C^2$, or continuous twice differentiable surface, almost everywhere \cite{stam1998evaluation,loop1987smooth}. It is this smoothness that has generated extensive interest in developing methods to describe physics \emph{although}, we note, NURBS is the industry standard for engineering design simply because it predates subdivision. Advances in subdivision based analysis methods have been numerous, ranging from analysis of thin shells \cite{green2002subdivision,liu2018isogeometric,cirak2001fully}, topology and shape optimization \cite{chen2022bi,chen2020acoustic,Jan_multiresolution,FEM_multires,alsnayyan2021laplacebeltrami}, design of subdivision exterior calculus \cite{de2016subdivision}, electromagnetic analysis \cite{Jie,alsnayyan2021isogeometric,Fu2017GeneralizedDS}, and acoustic analysis \cite{abdel_acoustics,chen2020subdivision}.

Our work on IGA on subdivision surfaces has taken two different paths: (a) defining currents via a Helmholtz decomposition on simply connected surfaces \cite{Jie,alsnayyan2021isogeometric} and (b) defining a div-conforming basis sets \cite{multiply_subd}. Using the former, we have been able to demonstrate numerous benefits of smoothness and thereby a complete Helmholtz decomposition
on simply connected surfaces. This permits seamless implementation of Calder\'{o}n operators \cite{Jie,alsnayyan2021isogeometric}, Debye sources \cite{Fu2017GeneralizedDS,epstein2010debye} and scalar integral equations \cite{li2018formulation}. They show excellent promise due to their accuracy. Exact Helmholtz decomposition on simply connected surfaces together with Calder\'{o}n-complex regularized combined field integral equation (CC-CFIER) enables us to easily overcome critical bottlenecks, specifically:  (i) ill-conditioning at low frequency; (ii) ill-conditioning due to mesh density; (iii) presence of spurious resonances; (iv) challenges due to multi-scale meshes \cite{Jie,alsnayyan2021isogeometric}. More recently, we have been able to solve a fundamental problem that plagues all higher order methods--cost of evaluation of near-field matrix elements \cite{alsnayyan2021isogeometric}. This cost is specially cumbersome when both the representation of the surface and currents is higher order. Additionally, we introduced the notion of manifold harmonics that is akin to a Fourier representation of quantities on the surface, and enables compression, and analysis on hierarchy of geometries \cite{alsnayyan2021isogeometric,alsnayyan2021laplacebeltrami}. 

What is critically missing from this suite is analysis of multiply connected structures. The harmonic basis functions or global loops on multiply connected structures are not easily handled within the purview of defining basis sets via surface derivatives, as they are both divergence-free and curl-free. One alternative is to avoid the need to define these decompositions by using divergence conforming basis \cite{multiply_subd}. But, we no longer will be able to take advantage of the properties described in the previous paragraph. Alternatively, we could use divergence conforming functions together with exact Helmholtz decomposition in a manner akin to that prescribed in  \cite{andriulli_precond} (i.e., use projectors), but at the expense of the number of degrees of freedom. In this paper, we present an alternate method for defining a harmonic basis by using vector fields wherein each component is a randomized polynomial \cite{Mike_Harmonic}. We will show that (a) the representation is exact in that the basis sets are divergence-free, curl-free, and both divergence-free and curl-free, (b) they are orthogonal to each other, and (c) the representation captures analytical fields defined on a torus and (d) scattering data obtained using these basis agrees well with those obtained using classical Rao-Wilton-Glisson based solutions.

The rest of the paper is organized as follows; in Section \ref{sec:probStatement} we define the problem. Next, we discuss the representation of surface, currents on the surface and their properties in Section \ref{sec:representation}. In Section \ref{sec:solvers} and \ref{sec:widebandFMM} we briefly discuss the necessary steps to discretize the system as well as a wideband MLFMA technique to rapidly evaluate the requisite inner products. Next, we present a number of results in Section \ref{sec:results} that validate the proposed approach. Finally, we summarize our findings in Section \ref{sec:summary} and outline future research avenues. 

\section{Problem statement\label{sec:probStatement}}

We consider the analysis of scattered fields $\{ \textbf{E}^{s},\textbf{H}^{s}\}$, from a perfect electrically conducting (PEC) object
$\Omega$, due to fields $\{ \textbf{E}^{i},\textbf{H}^{i}\}$ incident on its boundary $\Gamma \in \Omega$. It is assumed that this surface is equipped with a unique outward pointing normal denoted by $\hat{\textbf{n}}(\textbf{r})$, $\textbf{r} \in \Gamma$. The region external to this volume $\{\mathbb{R}^{3}\setminus \Omega\}$ is occupied by free space. The scattered field by the object at $\textbf{r} \in \{\mathbb{R}^{3}\setminus\Omega\}$ may be obtained using a Calder\'{o}n-Complex combined field integral equation (CC-CFIER) formulated in terms of the unknown surface current density $\vb{J}(\vb{r})$ on $\Gamma$ as follows:

\begin{equation}
    \label{eq:CCCFIER}
     \begin{split}
    & \left  (\frac{\mathcal{I}}{2} - \mathcal{K}_{\kappa} \right) \circ \textbf{J} (\vb{r})  - 2\mathcal{T}_{\kappa'} \circ \mathcal{T}_{\kappa} \circ \mathbf{J}(\vb{r}) =\\ & \hat{\textbf{n}} (\vb{r}) \cross \textbf{H}^{i}(\vb{r}) + 2\mathcal{T}_{\kappa'} \circ (\hat{\textbf{n}}(\vb{r}) \cross \textbf{E}^{i})(\vb{r}),
    \end{split}
\end{equation}
where,
\begin{subequations}
\begin{align}
\begin{split}
         \mathcal{T}_{\kappa} \circ \textbf{J}(\textbf{r}) &=  -j\eta \kappa\hat{\textbf{n}}(\textbf{r}) \times \int_{\Gamma} G_{\kappa}(\textbf{r},\textbf{r}')\cdot \textbf{J}(\textbf{r}')d\textbf{r}' \\  &+ j\frac{\eta}{\kappa}\hat{\textbf{n}}(\textbf{r}) \times \nabla\int_{\Gamma}  G_{\kappa}(\textbf{r},\textbf{r}') \nabla' \cdot \textbf{J}(\textbf{r}')d\textbf{r}',
\end{split}
\end{align}
\begin{equation}
    \mathcal{K}_{\kappa} \circ \textbf{J}(\textbf{r}) =  \hat{\textbf{n}}(\textbf{r}) \times \dashint_{\Gamma} \nabla G_{\kappa}(\textbf{r},\textbf{r}')\cdot \textbf{J}(\textbf{r}')d\textbf{r}'. 
\end{equation}
\end{subequations}
Here, $\mathcal{I}$ is the idempotent, $G_{\kappa}(\textbf{r},\textbf{r}') = \mbox{exp}[-j\kappa\abs{\textbf{r}-\textbf{r}'}]/(4 \pi \abs{\textbf{r}-\textbf{r}'})$, $\kappa$ is the free space wavenumber, $\eta$ is the free space impedance, $\kappa' = \kappa - j0.4\varsigma^{2/3}\kappa^{1/3}$ and $\varsigma$ is the maximum of the absolute values of mean curvatures on surface $\Gamma$, and $\mathcal{K}_{\kappa}$ is taken in the Cauchy principal value sense. In the above expressions, and what follows, we assume and suppress $\mbox{exp}[j\omega t]$ time dependence.

The construction of this formulation is derived from regularizing operators based on Calder\'{o}n identities and complexification techniques such that boundary integral operators on the left hand side of (\ref{eq:CCCFIER}) are second kind Fredholm operators. The construction of such regularizing operators has been proposed and analyzed in the literature \cite{boubendir2014well,DARBAS2006834}.

To solve \eqref{eq:CCCFIER} we will (i) represent the multiply-connected surface of the scatterer using isogeometric Loop subdivision basis sets, (ii) represent the currents on the surface using the \emph{same} basis set, and (iii) validate solutions to these integral equations solved using this procedure. Next, we discuss these in sequence.

\section{Representation of surface and currents\label{sec:representation}}   

In what follows, we will define subdivision surfaces, followed by definition of basis sets on these surfaces. Note, since we are developing an iso-geometric method, the same basis set will be used to define \emph{both} the manifold and the physics on the manifold. 

\subsection{Surface representation}

Here we provide a terse summary of Loop subdivision; information provided is purely for completeness and omits details that can be found in \cite{Ciraksubd,loop1987smooth} and references therein.  Let $T^{k}$ denote a $k$-th refined control mesh, with vertices $V^{k} := \{\textbf{v}_{i}, i = 1,\ldots, N_{v}\}$ and triangular faces $P^{k} := \{\textbf{p}_{i}, i = 1,\ldots, N_{f}\}$. Using this control mesh and approximating loop subdivision rules \cite{subdivision}, one can define a limit surface, $\Gamma (\vb{r})$, that is $C^2$ smooth almost everywhere, other than at some isolated points, using \emph{any} $k$-th refined control mesh $T^{k}$. Note the following: (a) a limit surface is obtained by successive refinement of a control mesh such that $k \longrightarrow \infty$, (b) there exists an analytical map from a control mesh to the limit surface, and (c) all control mesh refinement point to the same limit surface. We note that the construction of the limit surface is different from the conventional Lagrangian description. 

\begin{figure}[!t]
  \centering
  \includegraphics[width=3cm]{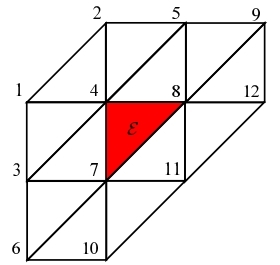}
  \caption{Regular triangular patch defined by its 1-ring vertices.
  \label{fig:1ring}}  
\end{figure}

Consider a patch/triangle $\epsilon$ and its set of control nodes, as shown in Fig. \ref{fig:1ring}. We define the following: (a) its 0-ring are the vertices that belong to the patch, (b) its 1-ring as the set of all vertices, $n_{v}$, that can be reached by traversing no more than two edges, and (c) the valence of a vertex is the number of edges that are incident on the vertex. In a typical Lagrangian description, the 0-th order description would be defined by the plane formed by the 0-ring. A higher order description would be due additional interpolation points within the 0-ring \cite{peterson1998computational}; this is in contrast to Loop subdivision. The map between a patch $\epsilon$ and the limit surface is defined by a weighted average of quartic box splines defined at its 1-ring neighbor vertices. This definition changes for vertices whose valence is \emph{not} 6, or irregular/extraordinary vertices \cite{stam1998evaluation}. Nevertheless, using this definition, on can define an effective basis function such that the entire surface $\Gamma (\vb{r})$ can be represented using \cite{stam1998evaluation,JieDaultShankerChapterSubd}
\begin{equation}
    \label{eq:Sexp}
        \Gamma(\textbf{r}) = \sum_{i=1}^{N_v} \textbf{c}_{i}\xi_{i}(\textbf{r}),
\end{equation}
where $\mathbf{c}_i$ is the $i$-th control node and $\xi_{i}$ is the effective basis function that is associated with $\mathbf{c}_i$ and has support $\Gamma_i$; note, $\cup_i \Gamma_i = \Gamma$.  The basis functions $\xi_{i}$ span a IGA finite dimensional space $\Psi$ that is the subspace of the Sobolev space $H^{2}(\Gamma)$ \cite{IGA_JCP,Rumpf_IGA}. We highlight the following properties of $\xi_i (\vb{r})$ that are noteworthy for defining basis sets that can model electromagnetic analysis: (a) positive in the domain of support, (b) compact support, (c) partition of unity, (d) the function and its surface gradient go smoothly to zero at the boundary of the domain of support and (e) $C^{2}$ continuity almost everywhere. 

\subsection{Representation of currents}

Consider representing a scalar function on a genus $g$ manifold $\Gamma$. To use isogeometric basis sets, we assume that there exists a net of control function values associated with each control vertex. In a manner similar to what is used in \eqref{eq:Sexp}, any scalar function can be expressed in terms of the Loop subdivision basis set via 
\begin{equation}
    f(\mathbf{r})= \sum_{n=1}^{N_v} a_{n}\xi_{n}(\mathbf{r}),
    \label{eq:scalar_func}
\end{equation}
where  $N_{v}$ and $\xi_{i}(\mathbf{r})$ are defined earlier. To define a basis for an electric current we note that on the manifold the current can be represented using a Helmholtz decomposition as 
\begin{equation}\label{eq:currentDecomp}
    \vb{J}(\vb{r}) = \nabla_\Gamma \phi (\vb{r}) + \hat{\vb{n}} (\vb{r}) \times \nabla_\Gamma \psi (\vb{r}) + \bar{\boldsymbol{\omega}} (\vb{r}),
\end{equation}
where $\nabla_\Gamma$ is surface gradient, $\phi (\vb{r})$ and $\psi (\vb{r})$ are potentials on the surface that satisfy the zero mean constraint, and $\bar{\omega}(\vb{r})$ is the harmonic component. Using the basis functions defined in \eqref{eq:scalar_func},  we can rewrite the potentials on the limit surface as
\begin{equation}
\begin{split}
     \phi(\textbf{r}) &\approx \tilde{\phi}(\textbf{r}) = \sum_{n=1}^{N_{v}}a^{1}_{n}\xi_{n}(\textbf{r}), \\  
    \psi(\textbf{r}) &\approx \tilde{\psi}(\textbf{r}) = \sum_{n=1}^{N_{v}}a^{2}_{n}\xi_{n}(\textbf{r}).
\end{split} 
\label{eq:f_approx}
\end{equation}
where $a_n^k$ for $k \in \{1,2\}$ are the unknown coefficients. It follows that the curl-free and divergence-free components of the current can be represented using the approximation of the potentials, viz., 
\begin{subequations}
\label{eq:J_Loop}
\begin{equation}
    \begin{split}
    &\textbf{J}(\textbf{r}) \approx \textbf{J}_{N}(\textbf{r}) = \sum_{n=1}^{N_{v}}\left[ a^{1}_{n}\textbf{J}_{n}^{1}(\textbf{r}) + a^{2}_{n}\textbf{J}_{n}^{2}(\textbf{r})\right] + \bar{\boldsymbol\omega}(\vb{r}),
    \end{split}
\end{equation}
\begin{equation}
    \begin{split}
    \label{eq:Loop_HH}
    &\textbf{J}_{n}^{1}(\textbf{r}) = \nabla_{\Gamma}\xi_{n}(\textbf{r}),\\
    &\textbf{J}_{n}^{2}(\textbf{r}) = \hat{\textbf{n}}(\textbf{r}) \times \nabla_{\Gamma}\xi_{n}(\textbf{r}),\\
    &\bar{\boldsymbol\omega}(\vb{r})= \sum_{n=1}^{g}\left[ a^{3}_{n}\textbf{J}_{n}^{3}(\textbf{r}) + a^{4}_{n}\textbf{J}_{n}^{4}(\textbf{r})\right],
    \end{split}
\end{equation}
\end{subequations}
wherein we define $\vb{J}_{n}^{3} (\vb{r}) $ and $\vb{J}_{n}^{4} (\vb{r}) $ in later sections. Note that since the representation is constructed using conditions on currents that rely on derivatives of the potentials $\tilde{\phi}(\textbf{r})$ and $\tilde{\psi}(\textbf{r})$, this leads to the existence of nontrivial solutions to (\ref{eq:J_Loop}) and therefore we must enforce uniqueness. In order to ensure uniqueness, we impose a zero-mean constraint for both divergence-free and curl-free components. A more thorough explanation, as well as, several properties of the basis functions can be found in \cite{Jie,Fu2017GeneralizedDS}.    

\subsection{Representation of harmonic components\label{sec:repHarm}}

Next, we detail the representation of the harmonic component. We begin by highlighting particular features of the Loop subdivision representation of the current thus far; more details can be found in \cite{Jie,JieDaultShankerChapterSubd}. They are as follows: 
\begin{enumerate}
    \item The representation of currents is $C^1$ continuous on $\Gamma (\vb{r})$ almost everywhere. 
    \item The normal is  $C^1$ almost everywhere. 
    \item On $\Gamma_n \cap \Gamma_m$, $\langle \vb{J}_n^1(\vb{r}), \vb{J}_m^2 (\vb{r})\rangle = -\langle \vb{J}_m^1(\vb{r}), \vb{J}_n^2 (\vb{r})\rangle$. 
    \item The inner product, $\langle \vb{J}_m^1 (\vb{r}), \vb{J}_n^2 (\vb{r})\rangle$, is identically zero. This can be proven using Green's theorems, $\xi_{m}(\vb{r})=0$ for $\vb{r} \in \partial\Gamma_{m}$, and properties stated earlier. To wit,
    \begin{equation}
        \begin{split}
             \langle \vb{J}_m^1 (\vb{r}), \vb{J}_n^2 (\vb{r})\rangle &= \int_{\Gamma_m \cap \Gamma_n} d\vb{r} \nabla_\Gamma \xi_m (\vb{r}) \cdot \vb{J}_n^2 (\vb{r})  \\
           &  = \int_{\partial(\Gamma_m \cap \Gamma_n)} \xi_m (\vb{r}) \hat{\vb{u}}(\vb{r})\cdot\vb{J}_n^2 (\vb{r})d\vb{r}  \\
            &- \int_{\Gamma_n \cap \Gamma_m} \xi_m (\vb{r}) \nabla_\Gamma \cdot \vb{J}_n^2 (\vb{r})d\vb{r} = 0,  
        \end{split}
    \end{equation}
    where $\partial \Gamma_m \cap \Gamma_n$ is the boundary of $\Gamma_m \cap \Gamma_n$, and $\hat{\vb{u}}(\vb{r})$ is a outward pointing normal tangential to the boundary. Given the properties stated earlier, the integrals are all zero. 
    \item As a result, the basis functions defined earlier are exaclty orthogonal. 
    \item Finally, we denote $\langle \vb{J}_m^1 (\vb{r}), \vb{J}_n^1 (\vb{r})\rangle = \langle \vb{J}_m^2 (\vb{r}), \vb{J}_n^2 (\vb{r})\rangle = \gamma_{mn} \delta_{lk}$ where, 
    \begin{equation}
        \gamma_{mn} = \int_{\Gamma_{n} \cap \Gamma_{n}} \nabla_{s}\xi_{n}(\vb{r}) \cdot \nabla_{s}\xi_{m}(\vb{r}) d\vb{r}. 
    \end{equation}
\end{enumerate}
These properties can be exploited readily in constructing a purely numerical basis for the harmonic components. We build on the ideas presented in \cite{Mike_Harmonic}. Prior to proceeding, we note that the space of harmonic fields is $2g$-dimensional, where $g$ is the genus of the object. Consider a function $\vb{F}(\vb{r})=-\hat{\vb{n}}(\vb{r}) \times \hat{\vb{n}}(\vb{r}) \times \vb{V}(\vb{r})$ where $\vb{V} (\vb{r}) \in \mathbb{R}^{3}$ and whose components are a low degree random polynomial. Specifically, let $\Gamma (\vb{r}) \in \left \{ (0,L_x) \times (0, L_y) \times (0, L_z) \right \} $. Then each component of $\vb{V}(\vb{r})$ can be defined as a tensor product of Legendre polynomials via $P_{\alpha}(x/2L_x)P_{\beta}(y/2L_y)P_{\zeta}(z/2L_z)$ where $P_{\zeta} (\cdot) $ is the standard Legendre polynomials of degree $\zeta$ defined on $[-1, 1]$.  Moreover, the maximum degree of Legendre polynomials should be at least $\ceil{(2g)^{1/3}}$ so that dimension of the space of polynomials is greater than the dimension of the space of the vector fields. For the examples in this paper, we choose the degree to be $\alpha + \beta + \zeta = \ceil{(2g)^{1/3}}+5$.

To obtain $g$ harmonic basis functions that project to $\bar{\boldsymbol{\omega}} (\vb{r})$, we start with the following: consider a collection of $g$ random vector fields formed using Legendre polynomials defined earlier and denoted using $\vb{F}_i (\vb{r})$ for $i \in \{1, \cdots, g\}$. For each such field $\vb{F}_i (\vb{r})$, we can obtain 
\begin{equation}\label{eq:harmBasis}
    \vb{J}^{3}_i (\vb{r}) = \vb{F}_i(\vb{r}) - \sum_n \left [ a_{n,i}^1\vb{J}^1_n (\vb{r})  + a_{n,i}^2\vb{J}^2_n (\vb{r}) \right ] 
\end{equation}
where the coefficients are arranged as $\mathcal{I}^1 = \left [a^1_{1,i},\cdots, a_{n,i}^1 \right ] ^T$ and $\mathcal{I}^2 = \left [a^2_{1,i}, \cdots, a_{n,i}^2 \right ]^T$, and obtained by solving 
\begin{subequations}\label{eq:helmCompnents}
\begin{equation}
    \label{eq:Gram_system}
    \begin{bmatrix}
    \mathcal{G}^{11} & 0\\
    0 & \mathcal{G}^{22} \\
    \end{bmatrix}
    \begin{bmatrix}
    \mathcal{I}^{1}\\
    \mathcal{I}^{2} \\
    \end{bmatrix}
    =
    \begin{bmatrix}
    \mathcal{V}^{1}_i\\
    \mathcal{V}^{2}_i \\
    \end{bmatrix}.\\
\end{equation}
Here, $\mathcal{G}^{11}_{nm} = \mathcal{G}^{22}_{nm} = \gamma_{nm}$, and 
\begin{equation}
    \mathcal{V}^{k}_{n,i} = \int_{\Gamma_{n}} \vb{J}_{n}^{k}(\vb{r}) \cdot \vb{F}_i(\vb{r})d\vb{r}.    
\end{equation}
\end{subequations}
 Furthermore, note that if $\vb{J}^{3}_i (\vb{r})$ is a harmonic vector field, then $\vb{J}^{4}_i (\vb{r}) = \hat{\vb{n}}(\vb{r}) \times \vb{J}^{3}_i (\vb{r})$ is also a harmonic vector field that is independent of $\vb{J}^{3}_i (\vb{r})$. These basis span the space of harmonic fields (with high probability that the set of random vector fields are linearly independent) \cite{Mike_Harmonic}.

\subsection{More properties of basis sets\label{sec:repHarm2}}

There are several properties that make this decomposition viable for analysis and as such we examine these properties below. Specifically, we note the following: 
\begin{enumerate}
    \item The inner products of a field $\vb{F} (\vb{r})$ with basis $\vb{J}_{n}^{1}(\vb{r})$ and $\vb{J}_{n}^2(\vb{r})$, for $n \in \{1,\ldots,N_{v}\}$, effectively measures the surface divergence and curl of the field.
    \item It it trivial to show that both surface divergence and surface curl of $\vb{J}_{n}^{3} (\vb{r})$ and $\vb{J}_{n}^{4}(\vb{r})$, for $n \in \{1,\ldots,g\}$, as defined in \eqref{eq:harmBasis} is identically zero via \eqref{eq:helmCompnents}. 
    \item Along the same lines, it can be shown that $\langle \vb{J}_n^k (\vb{r}), \vb{J}^{l}_i (\vb{r}) \rangle = 0$ for $k \in \{1, 2\}$ and $l \in \{3, 4\}$.
    \item The support of $\vb{J}^{3}_{i}(\vb{r})$ and  $\vb{J}^{4}_{i}(\vb{r})$ is global and as such we need additional infrastructure to evaluate the necessary matrix elements. 
\end{enumerate}
A direct consequence of these properties and those in section \ref{sec:repHarm} is that the prescribed basis $\vb{J}_n^k(\vb{r})$ for $n \in \{1, \dots, N_v\}$ and $k \in \{1,2\}$, and $\vb{J}^{l}_i (\vb{r})$ for $l \in \{3,4\}$ and $i \in \{1, \cdots, g\}$, provide a complete Helmholtz decomposition of vector quantities defined on the manifold. 

\subsection{Fourier representation of and on the manifold}
This subsection is not meant to be exhaustive, but to point to other potential applications that tie in into the topics of current interest. We begin by remarking that the  Laplace–Beltrami operator \cite{Italian_LBO_book} admits a complete and countable sequence of eigenvectors which form an orthonormal basis in $L_{2}\left(\Gamma\right)$, otherwise known as manifold harmonics; MHB is tantamount to Fourier basis on the manifold. As such, these manifold harmonics are used extensively in computer graphics for their spectral representation of shapes (even complex ones), as it provides a compressed representation and is robust to deformation \cite{Italian_LBO_paper}. Given that the manifold harmonics can be used to define potentials $\phi(\vb{r})$ and $\psi(\vb{r})$ on the surface, we have used manifold representation for both analysis \cite{alsnayyan2021isogeometric} and shape optimization \cite{alsnayyan2021laplacebeltrami}. It stands to reason that as the proposed technique leverages representation of $\phi (\vb{r}) $ and $\psi (\vb{r})$ one can readily extend the manifold harmonic based approach to analysis of multiply connected objects. In addition, this technique can be useful in various computer graphics applications including spectral processing of tangential vector fields \cite{azencot2013operator,spectral_tangential,poelke2017hodge,custers2020subdivision}.

\section{Field Solvers\label{sec:solvers}}

Thus far, we have discussed the representation of current on the surface. We note that all basis function are defined for all values of $\vb{r} \in \Gamma (\vb{r})$ and the harmonic basis have global support. We detail the discretization of \eqref{eq:CCCFIER}, in terms of the prescribed basis set; in particular, we use a Galerkin scheme to discretize these equations. We note that efficient evaluation of inner products is discussed in the next section. Note that discretizing Calder\'{o}n type operators requires intermediate spaces, effected through a Gram matrix. We define the required Gram-matrix $[G]$ using
\begin{subequations}
\begin{equation}
[G] = 
\begin{bmatrix}
\mathcal{G}^{11} & 0 & 0&0\\
0 & \mathcal{G}^{22} & 0&0\\
 0& 0 & \mathcal{G}^{33}&\mathcal{G}^{34}\\
 0& 0 & \mathcal{G}^{43}&\mathcal{G}^{44}
\end{bmatrix},
\end{equation}
\begin{equation}
    G^{lk}_{nm} = \langle \vb{J}^{l}_n (\vb{r}), \vb{J}^{k}_m(\vb{r}) \rangle = \int_{\Gamma} \vb{J}^{l}_n (\vb{r}) \cdot \vb{J}^{k}_m(\vb{r}) d\vb{r},
\end{equation}
\end{subequations}
where $l,k \in \{3,4\}$ and $G^{33}_{nm} = G^{44}_{nm}$, $G^{34}_{nm} = - G^{43}_{mn} $. Using \eqref{eq:J_Loop} in \eqref{eq:CCCFIER} and Galerkin testing results in a matrix system that can be written as 
\begin{subequations}
\begin{equation}
    \label{eq:matrix_eq}
    \left [ Z\right ] \left [I \right ] = \left [ V\right ]
\end{equation}
where,
\begin{equation}
    \left [ Z \right] =[G]^{-1}\left [  \left [ L \right ] + \left [ K \right ] \right ]
\end{equation}
with
\begin{equation}
[K]_{nm}^{lk} = \left \langle \textbf{J}_n^l (\vb{r}), \frac{\textbf{J}^k_m}{2} (\textbf{r}) - \mathcal{K}_{\kappa}\circ \textbf{J}_m^k(\textbf{r}) \right \rangle_{\Gamma_{n}},
\end{equation}
\begin{equation}
    \label{eq:T}
    [T]_{\widetilde{\kappa},nm}^{lk} = \left \langle \textbf{J}_n^l (\vb{r}),  \mathcal{T}_{\widetilde{\kappa}}\circ \textbf{J}_m^k(\textbf{r}) \right \rangle_{\Gamma_{n}},
\end{equation}
where $\widetilde{\kappa} \in \{\kappa',\kappa\}$, and, as defined earlier $ \kappa' = \kappa + 0.4\varsigma^{2/3}\kappa^{1/3}$, and $\varsigma$ is the maximum of the magnitude of the mean curvature of the object, and
\begin{equation}
[L] = -2[T]_{\kappa'}[G]^{-1}[T]_\kappa.
\end{equation}
\end{subequations}
Furthermore, we have
\begin{subequations}
\begin{equation}
    [I]_{m}^{k} = a^k_{m},
\end{equation}

\begin{equation}
    [V]_{n}^{k} = [G]^{-1}\left[-2[T]_{\kappa'}[G]^{-1}[V_{T}]_{n}^{k} +[V_{K}]_{n}^{k}\right],  
\end{equation}
\end{subequations}
with  
\begin{subequations}
\begin{equation}
    [V_{T}]_{n}^{k}  = \left \langle \textbf{J}_{n}^{k}(\textbf{r}),   \textbf{E}^{i}(\textbf{r})\right \rangle_{\Gamma_{n}},
\end{equation}
\begin{equation}
    [V_{K}]_{n}^{k} = \left \langle \textbf{J}_{n}^{k}(\textbf{r}),\hat{\textbf{n}}  \cross \textbf{H}^{i}(\textbf{r})\right \rangle_{\Gamma_{n}},  
\end{equation}
\end{subequations}
Here, we have defined $l,k \in \{1,2,3,4\}$. Lastly, we note that the stabilizing properties of the Calderón preconditioner are local \cite{Local_CC_CFIER}, which allows the use of a localized version of the preconditioner $[T]_{\kappa'}$. As such, we choose to omit all interactions of a distance greater than $1.25\lambda$.  

\section{Wideband MLFMA for Evaluation of Inner Products \label{sec:widebandFMM}}

One of the critical components is the evaluation of inner products between different types of basis functions. The principal challenge arises from the higher order nature of the basis set as well as the domain of support. Indeed, the size of the domain $\Gamma_n$ associated with a basis function can be as large as $0.9\lambda$, if not larger. This is a direct cause of the fourth-order geometry and the basis set being third order. In addition, the harmonic basis has global support. As a consequence of these attributes, one needs a specialized fast method equipped to effectively evaluate all required inner products, given their higher order quadrature rules. In this paper, we leverage a mixed potential approach together with wideband MLFMA. This approach exploits the structure of wideband MLFMA coupled with an adaptive integration rule that relies on partitioning each patch into sub-triangles. Specifically, consider a basis function interacting with itself; the inner product devolves to a collection of sub-patches, such that each sub-patch corresponds to a low order integration rule and sub-patch interactions can be partitioned into near and far regions via a tree based algorithm, with leaf box sizes as small as the smallest sub-patch, i.e. $0.0625\lambda$ or smaller \cite{vikram2009novel,dault2016mixed}. This approach has been developed for higher order methods \cite{dault2016mixed} and applied to subdivision basis sets for both acoustic and electromagnetic integral equations \cite{abdel_acoustics,alsnayyan2021isogeometric}. As is evident from the results in these paper, we have an error controllable method for evaluating all matrix vector products (both in the near and far field of each other). 

\section{Numerical Results\label{sec:results}}    

In this section, numerical examples are given to demonstrate the convergence and effectiveness of the presented numerical scheme to electromagnetic analysis. In order to do so, we shall present data on the following: (i) the exact Helmholtz decomposition of analytical fields defined on a torus; (ii) application of our wideband MLFMA scheme for analyzing complex multiply-connected structures; (iii) timings and iteration count for overall GMRES iterative solve. 

\begin{figure}
    \centering
    \includegraphics[width=0.5\textwidth]{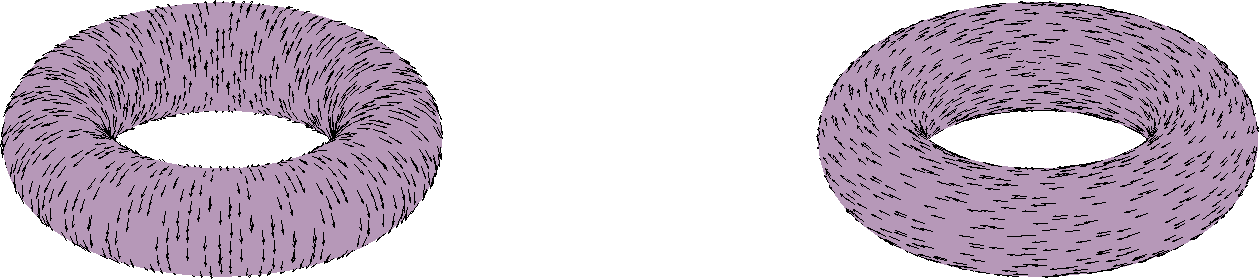}
    \caption{ The two basis harmonic vector fields for a torus.}
    \label{fig:harmonic_fields}
\end{figure}    

\subsection{Helmholtz Decomposition}

In our first example, we demonstrate the effectiveness of our prescribed basis set to affect a Helmholtz decomposition of a vector field that is a linear combination of analytically prescribed curl-free, divergence-free, and harmonic vector fields, defined on a torus. We choose a torus, see Fig.~\ref{fig:harmonic_fields}, with a major axis of 3 m and a minor axis of 1 m, modeled using an initial control mesh comprising of 2048 vertices and 4096 faces. We refine the control mesh twice such that surface area and the total curvature of the limit surface $\Gamma$ agree within 99\% to the analytical torus. We aim to report the $L_{2}$ norm of each component of the Helmholtz decomposition, for the following vector field $\bar{\vb{X}}(\vb{r})$. In particular, we have
\begin{align}\label{eq:testHelm}
\begin{split}
        \vb{X}_{d}(\vb{r}) &= -\hat{\vb{n}}(\vb{r}) \times \hat{\vb{n}}(\vb{r}) \times  \vb{r}, \\
            \vb{X}_{c}(\vb{r}) &= \hat{\vb{n}} \times \vb{X}_{d}(\vb{r}),\\
               \textbf{X}_{h}(\vb{r}) &= -\hat{\vb{n}}(\vb{r}) \times \hat{\vb{n}}(\vb{r}) \times \frac{\vb{t}}{||\vb{t}||^{2}},\\
               \bar{\vb{X}}(\vb{r}) &=
               \vb{X}_{d}(\vb{r}) + \vb{X}_{c}(\vb{r}) + \textbf{X}_{h}(\vb{r}) + \hat{\vb{n}}(\vb{r})  \times  \textbf{X}_{h}(\vb{r}) 
\end{split}
\end{align} 
where $\vb{t} = y \hat{x} - x \hat{y}$, and $\vb{X}_{c}$ is a divergence-free field, $\vb{X}_{d}$ is a curl-free field, and $\vb{X}_{h}$ and $\hat{\vb{n}}\times\vb{X}_h$ are purely harmonic fields. Next, we choose to approximate these fields using the basis sets defined earlier; specifically, 
\begin{align}
\begin{split}
        \bar{\vb{X}}(\vb{r}) &\approx \vb{J}^{1}(\vb{r}) + \vb{J}^{2}(\vb{r}) + \vb{J}^{3}(\vb{r}) + \vb{J}^{4}(\vb{r}) \\
            &= \sum_{n} a_{n}^{1} \nabla_{\Gamma}\vb{J}^{1}_{n}(\vb{r}) + \sum_{n}a_{n}^{2} \vb{J}^{2}_{n}\xi_{n}(\vb{r}) \\
            &+ \sum_{n=1}^{g}a^{3}_{n}\vb{J}^{3}_{n}(\vb{r})
            + \sum_{n=1}^{g}a^{4}_{n}\vb{J}^{4}_{n}(\vb{r}),
\end{split}
\end{align}    
where, $\vb{J}^{1}(\vb{r}) \approx \vb{X}_{d}(\vb{r})$, $\vb{J}^{2}(\vb{r}) \approx \vb{X}_{c}(\vb{r})$, $\vb{J}^{3}(\vb{r}) \approx \vb{X}_{h}(\vb{r})$, and $\vb{J}^{4}(\vb{r}) \approx \hat{\vb{n}}(\vb{r}) \times \vb{X}_{h}(\vb{r})$. We compute the unknown coefficients using the orthogonality properties of the basis set as presented in sections \ref{sec:repHarm} and \ref{sec:repHarm2}. 

\begin{table}[!h]
\begin{center}

\resizebox{1.0\columnwidth}{!}{
\begin{tabular}{c | ccccc}
 & $\vb{X}_{d}$ & $\vb{X}_{c}$ &  $\vb{X}_{h}$ &  $\hat{\vb{n}} \times \vb{X}_{h}$  \\
 \hline
 $\parallel \vb{J}^{1} \parallel_{2}$ & 5.05E2 & 2.28E-17 & 3.02E-10 & 3.31E-19 &  \\
 $\parallel\vb{J}^{2}\parallel_{2}$  &2.28E-17 & 5.05E2 & 3.31E-19 & 3.02E-10 & \\
 $\parallel\vb{J}^{3}\parallel_{2}$  &  1.01E-26 &1.67E-26 & 1.36E2 & 5.26E-18 &\\
 $\parallel\vb{J}^{4}\parallel_{2}$  & 1.98E-29 & 6.87E-24 & 5.02E-17 & 1.36E2 &
\end{tabular}
}
\end{center}
\caption{$L_{2}$-norms of the components of the decomposed vector fields}
\label{table1}
\end{table}

Table.~\ref{table1} lists the $L_{2}$-norm of each component of the decomposition of the vector fields, respectively. Furthermore, in Fig.~\ref{fig:harmonic_fields} we plot the two resulting harmonic fields. Ideally, the four components are recovered in the corresponding discrete spaces. As expected, the dominant components in each decomposition are reflected in the correct subspaces. The projection on to other subspaces are significantly smaller in magnitude, but do not vanish. This is largely due the fact that while \eqref{eq:testHelm} is analytically exact on a canonical torus, it is not so on a geometric approximation.

\subsection{Scattering from Complex Objects}
In what follows, we analyze the performance of two different cases of the proposed approach for studying the radar cross section (RCS) and compare it with a conventional MoM CFIE solution technique that relies on RWG basis functions, otherwise referred to as RWG-CFIE. These two cases are (a) the prescribed formulation in this paper, which we refer to as CC-CFIER: Loop + GL and (b) the same prescribed formulation, but we omit the harmonic component of the surface current, which we refer to as CC-CFIER: Loop; in other words we set $\vb{J}^{k}_{i}=0$, for $k \in \{3,4\}$ and $i\in\{1,\ldots,g\}$. In particular, we compare all three generated RCSs against each other, we report the iteration count to reach the specified GMRES solver tolerance and the time taken to reach the prescribed tolerance. Unless otherwise stated, we are comparing RCS in the $\phi = 0$ plane, due to a  plane wave field propagating in $\hat{\kappa} = -\hat{z}$ and polarized along $\hat{x}$ axis. Furthermore, in the experiments discussed next, the finer discretization was used with RWG basis (together with a Lagrangian geometry description). We ensured that the surface areas of the Lagrangian mesh agree within 99\% to the subdivision limit surface. All of the numerical results presented in the graphs in this section were obtained by prescribing a GMRES residual tolerance equal to $10^{-7}$ for the overall system and $10^{-11}$ for inverting the gram matrix.

\begin{figure}[!t]
    \centering
    \includegraphics[width=0.20\textwidth]{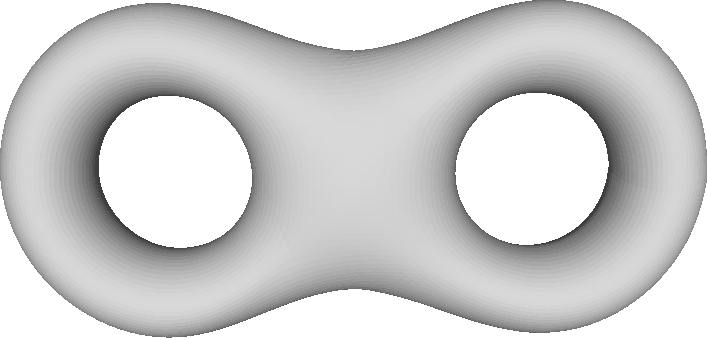}
     \caption{Double torus with 2 global loops.}
    \label{fig:Double_torus}
\end{figure}  

\begin{figure}[!h]
    \centering
    \includegraphics[width=0.4\textwidth]{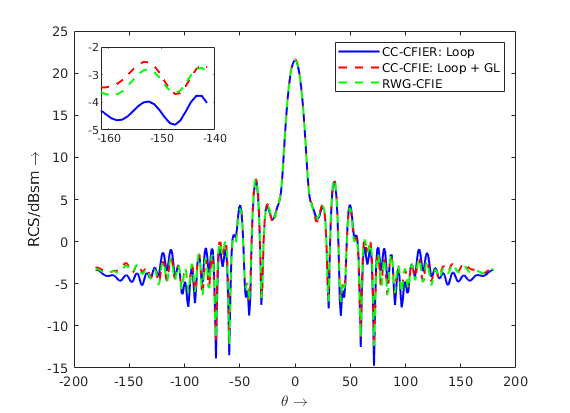}
    \caption{Radar cross section of the double torus ($\phi$ = 0 cut)}
    \label{fig:double_torus_rcs}
\end{figure}

The first example we consider is a double torus of genus 2, see Fig.~\ref{fig:Double_torus}, that fits in a $10\lambda \times 4.81\lambda \times 2.06\lambda$ box. The number of DoF for the CFIE: RWG is 36864, CC-CFIE: Loop + GL results in 6146 DoF, and lastly the CC-CFIE: Loop has 6142 DoF, respectively. The CFIE: RWG converges in 77 iterations after 2 m 7 s, CC-CFIE: Loop + GL converges in 17 iterations after 1 m 8 s, and CC-CFIE: Loop reaches tolerance within 22 iterations, for a total of 1 m 7 s. From Fig.~\ref{fig:double_torus_rcs}, we report excellent agreement between both the CC-CFIE: Loop + GL and the RWG-CFIE and we have highlighted the glaring disparity when omitting the harmonic component in CC-CFIE: Loop between the other methods in the inset.

\begin{figure}[!h]
    \centering
    \includegraphics[width=0.20\textwidth]{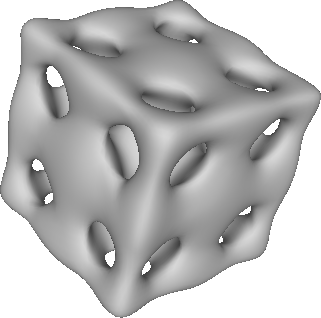}
     \caption{Chmutov surface with 28 global loops.}
    \label{fig:Chmutov}
\end{figure}

\begin{figure}[!h]
    \centering
    \includegraphics[width=0.4\textwidth]{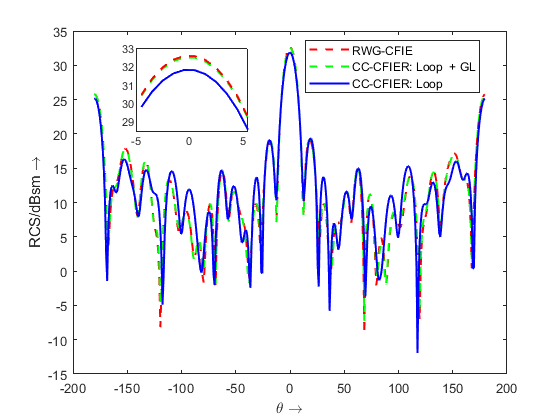}
     \caption{Radar cross section of the Chmutov ($\phi$ = 0 cut)}
     \label{fig:Chmutov_rcs}
\end{figure}

The second example is a Chmutov surface, see Fig.~\ref{fig:Chmutov}, that fits in a 5$\lambda$ box with a genus of 28. The number of DoF for the CFIE: RWG is 144000, CC-CFIE: Loop + GL has 23950 DoF, and lastly the CC-CFIE: Loop results in 23894 DoF, respectively. The CFIE: RWG converges after 350 iterations in 26 m 25s, CC-CFIE: Loop + GL converges after 36 iterations in 3 m 48 s, and lastly, the CC-CFIE: Loop converges after 104 iterations in 10 m 20 s. Fig.~\ref{fig:Chmutov_rcs} shows excellent agreement between both the CC-CFIE: Loop + GL and the RWG-CFIE and lastly, we highlight the glaring disparity when omitting the harmonic component in CC-CFIE: Loop between the other methods in the inset.

\section{Summary\label{sec:summary}}

This work develops the numerical infrastructure on to existing Loop subdivision IGA methods to enable analysis of multiply-connected structures. This is effected with the help of random vector projections onto subdivision surfaces. As a result, we are able to prescribe a complete Helmholtz decomposition on the surface wherein orthogonality of different components is preserved. This basis set can then be used to discretize well conditioned integral equations. A number of examples presented validate the assertions as well demonstrate the applicability to analysis of scattering from perfectly conducting objects. It is straightforward to apply the proposed basis set to other types of integral equations, such as decoupled potential/field integral equations. However, some problems remain that are a feature of the approximating subdivision basis we have chosen (and not subdivision in itself), i.e., preserving sharp features\cite{song2009modeling}. Our future goal is to enable the inclusion of divergence-conforming basis within this analysis rubric; the pathway is well known from our earlier work \cite{JieDaultShankerChapterSubd}. This will be presented in different forums soon.

\bibliographystyle{IEEEtran}
\bibliography{IEEEabrv,Bibliography}

\begin{thebibliography}{10}
\providecommand{\url}[1]{#1}
\csname url@rmstyle\endcsname
\providecommand{\newblock}{\relax}
\providecommand{\bibinfo}[2]{#2}
\providecommand\BIBentrySTDinterwordspacing{\spaceskip=0pt\relax}
\providecommand\BIBentryALTinterwordstretchfactor{4}
\providecommand\BIBentryALTinterwordspacing{\spaceskip=\fontdimen2\font plus
\BIBentryALTinterwordstretchfactor\fontdimen3\font minus
  \fontdimen4\font\relax}
\providecommand\BIBforeignlanguage[2]{{%
\expandafter\ifx\csname l@#1\endcsname\relax
\typeout{** WARNING: IEEEtran.bst: No hyphenation pattern has been}%
\typeout{** loaded for the language `#1'. Using the pattern for}%
\typeout{** the default language instead.}%
\else
\language=\csname l@#1\endcsname
\fi
#2}}
\renewcommand\BIBentryALTinterwordstretchfactor{4}

\bibitem{swanson2003microwave}
D.~G. Swanson and W.~J. Hoefer, \emph{Microwave circuit modeling using
  electromagnetic field simulation}.\hskip 1em plus 0.5em minus 0.4em\relax
  Artech House, 2003.

\bibitem{sideris2019ultrafast}
C.~Sideris, E.~Garza, and O.~P. Bruno, ``Ultrafast simulation and optimization
  of nanophotonic devices with integral equation methods,'' \emph{ACS
  Photonics}, vol.~6, no.~12, pp. 3233--3240, 2019.

\bibitem{garza2020boundary}
E.~Garza, ``Boundary integral equation methods for simulation and design of
  photonic devices,'' Ph.D. dissertation, 2020.

\bibitem{peterson1998computational}
A.~F. Peterson, S.~L. Ray, R.~Mittra, I.~of~Electrical, and E.~Engineers,
  \emph{Computational methods for electromagnetics}.\hskip 1em plus 0.5em minus
  0.4em\relax IEEE press New York, 1998, vol. 351.

\bibitem{vico2016decoupled}
F.~Vico, M.~Ferrando, L.~Greengard, and Z.~Gimbutas, ``The decoupled potential
  integral equation for time-harmonic electromagnetic scattering,''
  \emph{Communications on Pure and Applied Mathematics}, vol.~69, no.~4, pp.
  771--812, 2016.

\bibitem{liu2018potential}
Q.~S. Liu, S.~Sun, and W.~C. Chew, ``A potential-based integral equation method
  for low-frequency electromagnetic problems,'' \emph{IEEE transactions on
  antennas and Propagation}, vol.~66, no.~3, pp. 1413--1426, 2018.

\bibitem{hsiao1997mathematical}
G.~C. Hsiao and R.~E. Kleinman, ``Mathematical foundations for error estimation
  in numerical solutions of integral equations in electromagnetics,''
  \emph{IEEE transactions on Antennas and Propagation}, vol.~45, no.~3, pp.
  316--328, 1997.

\bibitem{calderon_EFIE}
S.~{Yan}, J.~{Jin}, and Z.~{Nie}, ``Efie analysis of low-frequency problems
  with loop-star decomposition and calderón multiplicative preconditioner,''
  \emph{IEEE Transactions on Antennas and Propagation}, vol.~58, no.~3, pp.
  857--867, 2010.

\bibitem{adrian2021electromagnetic}
S.~B. Adrian, A.~Dely, D.~Consoli, A.~Merlini, and F.~P. Andriulli,
  ``Electromagnetic integral equations: Insights in conditioning and
  preconditioning,'' \emph{IEEE Open Journal of Antennas and Propagation},
  2021.

\bibitem{boubendir2014well}
Y.~Boubendir and C.~Turc, ``Well-conditioned boundary integral equation
  formulations for the solution of high-frequency electromagnetic scattering
  problems,'' \emph{Computers \& Mathematics with Applications}, vol.~67,
  no.~10, pp. 1772--1805, 2014.

\bibitem{warnick2008numerical}
K.~F. Warnick, \emph{Numerical analysis for electromagnetic integral
  equations}.\hskip 1em plus 0.5em minus 0.4em\relax Artech House, 2008.

\bibitem{takahashi2022isogeometric}
T.~Takahashi, T.~Hirai, H.~Isakari, and T.~Matsumoto, ``An isogeometric
  boundary element method for three-dimensional doubly-periodic layered
  structures in electromagnetics,'' \emph{Engineering Analysis with Boundary
  Elements}, vol. 136, pp. 37--54, 2022.

\bibitem{wolf2021isogeometric}
F.~Wolf, ``Isogeometric boundary elements,'' in \emph{Analysis and
  Implementation of Isogeometric Boundary Elements for Electromagnetism}.\hskip
  1em plus 0.5em minus 0.4em\relax Springer, 2021, pp. 35--71.

\bibitem{Jie}
J.~Li, D.~Dault, B.~Liu, Y.~Tong, and B.~Shanker, ``Subdivision based
  isogeometric analysis technique for electric field integral equations for
  simply connected structures,'' \emph{Journal of computational physics}, vol.
  319, pp. 145--162, 2016.

\bibitem{simpson2018isogeometric}
R.~N. Simpson, Z.~Liu, R.~Vazquez, and J.~A. Evans, ``An isogeometric boundary
  element method for electromagnetic scattering with compatible b-spline
  discretizations,'' \emph{Journal of Computational Physics}, vol. 362, pp.
  264--289, 2018.

\bibitem{buffa2010isogeometric}
A.~Buffa, G.~Sangalli, and R.~V{\'a}zquez, ``Isogeometric analysis in
  electromagnetics: B-splines approximation,'' \emph{Computer Methods in
  Applied Mechanics and Engineering}, vol. 199, no. 17-20, pp. 1143--1152,
  2010.

\bibitem{dolz2019isogeometric}
J.~D\"olz, S.~Kurz, S.~Sch\"ops, and F.~Wolf, ``Isogeometric boundary elements
  in electromagnetism: rigorous analysis, fast methods, and examples,''
  \emph{SIAM Journal on Scientific Computing}, vol.~41, no.~5, pp. B983--B1010,
  2019.

\bibitem{Bazilevs2010}
Y.~Bazilevs, V.~M. Calo, J.~A. Cottrell, J.~A. Evans, T.~J.~R. Hughes,
  S.~Lipton, M.~A. Scott, and T.~W. Sederberg, ``Isogeometric analysis using
  t-splines,'' \emph{Computer methods in applied mechanics and engineering},
  vol. 199, no.~5, pp. 229--263, 2010.

\bibitem{zorin2007subdivision}
D.~Zorin, ``Subdivision on arbitrary meshes: algorithms and theory,'' in
  \emph{Mathematics and Computation in Imaging Science and Information
  Processing}.\hskip 1em plus 0.5em minus 0.4em\relax World Scientific, 2007,
  pp. 1--46.

\bibitem{stam1998evaluation}
J.~Stam, ``Evaluation of loop subdivision surfaces,'' \emph{SIGGRAPH 99 Course
  Notes}, 06 2001.

\bibitem{loop1987smooth}
C.~Loop, ``Smooth subdivision surfaces based on triangles,'' \emph{Master's
  thesis, University of Utah, Department of Mathematics}, 1987.

\bibitem{green2002subdivision}
S.~Green, G.~Turkiyyah, and D.~Storti, ``Subdivision-based multilevel methods
  for large scale engineering simulation of thin shells,'' in \emph{Proceedings
  of the seventh ACM symposium on Solid modeling and applications}, 2002, pp.
  265--272.

\bibitem{liu2018isogeometric}
Z.~Liu, M.~Majeed, F.~Cirak, and R.~N.~Simpson, ``Isogeometric fem-bem coupled
  structural-acoustic analysis of shells using subdivision surfaces,''
  \emph{International Journal for Numerical Methods in Engineering}, vol. 113,
  no.~9, pp. 1507--1530, 2018.

\bibitem{cirak2001fully}
F.~Cirak and M.~Ortiz, ``Fully c1-conforming subdivision elements for finite
  deformation thin-shell analysis,'' \emph{International Journal for Numerical
  Methods in Engineering}, vol.~51, no.~7, pp. 813--833, 2001.

\bibitem{chen2022bi}
L.~Chen, H.~Lian, Z.~Liu, Y.~Gong, C.~Zheng, and S.~Bordas, ``Bi-material
  topology optimization for fully coupled structural-acoustic systems with
  isogeometric fem--bem,'' \emph{Engineering Analysis with Boundary Elements},
  vol. 135, pp. 182--195, 2022.

\bibitem{chen2020acoustic}
L.~Chen, C.~Lu, H.~Lian, Z.~Liu, W.~Zhao, S.~Li, H.~Chen, and S.~P. Bordas,
  ``Acoustic topology optimization of sound absorbing materials directly from
  subdivision surfaces with isogeometric boundary element methods,''
  \emph{Computer Methods in Applied Mechanics and Engineering}, vol. 362, p.
  112806, 2020.

\bibitem{Jan_multiresolution}
J.~Zapletal and J.~Bouchala, ``Shape optimization and subdivision surface based
  approach to solving 3d bernoulli problems,'' \emph{Computers \& mathematics
  with applications (1987)}, vol.~78, no.~9, pp. 2911--2932, 2019.

\bibitem{FEM_multires}
K.~Bandara and F.~Cirak, ``Isogeometric shape optimisation of shell structures
  using multiresolution subdivision surfaces,'' \emph{Computer aided design},
  vol.~95, pp. 62--71, 2018.

\bibitem{alsnayyan2021laplacebeltrami}
A.~Alsnayyan and B.~Shanker, ``Laplace-beltrami based multi-resolution shape
  reconstruction on subdivision surfaces,'' \emph{arXiv preprint
  arXiv:2104.04027}, 2021.

\bibitem{de2016subdivision}
F.~De~Goes, M.~Desbrun, M.~Meyer, and T.~DeRose, ``Subdivision exterior
  calculus for geometry processing,'' \emph{ACM Transactions on Graphics
  (TOG)}, vol.~35, no.~4, pp. 1--11, 2016.

\bibitem{alsnayyan2021isogeometric}
A.~Alsnayyan and B.~Shanker, ``Iso-geometric integral equation solvers and
  their compression via manifold harmonics,'' \emph{arXiv preprint
  arXiv:2106.11907}, 2021.

\bibitem{Fu2017GeneralizedDS}
X.~Fu, J.~Li, L.~Jiang, and B.~Shanker, ``Generalized debye sources-based efie
  solver on subdivision surfaces,'' \emph{IEEE Transactions on Antennas and
  Propagation}, vol.~65, pp. 5376--5386, 2017.

\bibitem{abdel_acoustics}
A.~M.~A. Alsnayyan, J.~Li, S.~Hughey, A.~Diaz, and B.~Shanker, ``Efficient
  isogeometric boundary element method for analysis of acoustic scattering from
  rigid bodies,'' \emph{The Journal of the Acoustical Society of America}, vol.
  147, no.~5, pp. 3275--3284, 2020.

\bibitem{chen2020subdivision}
L.~Chen, C.~Lu, W.~Zhao, H.~Chen, and C.~Zheng, ``Subdivision
  surfaces—boundary element accelerated by fast multipole for the structural
  acoustic problem,'' \emph{Journal of Theoretical and Computational
  Acoustics}, vol.~28, no.~02, p. 2050011, 2020.

\bibitem{multiply_subd}
J.~Li and B.~Shanker, ``Isogeometric analysis of {EM} scattering on
  multiply-connected subdivision surfaces.''\hskip 1em plus 0.5em minus
  0.4em\relax 2017 IEEE International Symposium on Antennas and Propagation \&
  USNC/URSI National Radio Science Meeting, 2017, pp. 1557--1558.

\bibitem{epstein2010debye}
C.~L. Epstein and L.~Greengard, ``Debye sources and the numerical solution of
  the time harmonic maxwell equations,'' \emph{Communications on Pure and
  Applied Mathematics: A Journal Issued by the Courant Institute of
  Mathematical Sciences}, vol.~63, no.~4, pp. 413--463, 2010.

\bibitem{li2018formulation}
J.~Li, X.~Fu, and B.~Shanker, ``Formulation and iso-geometric analysis of
  scalar integral equations for electromagnetic scattering,'' \emph{IEEE
  Transactions on Antennas and Propagation}, vol.~66, no.~4, pp. 1957--1966,
  2018.

\bibitem{andriulli_precond}
F.~P. {Andriulli}, K.~{Cools}, H.~{Bagci}, F.~{Olyslager}, A.~{Buffa},
  S.~{Christiansen}, and E.~{Michielssen}, ``A multiplicative calderon
  preconditioner for the electric field integral equation,'' \emph{IEEE
  Transactions on Antennas and Propagation}, vol.~56, no.~8, pp. 2398--2412,
  2008.

\bibitem{Mike_Harmonic}
D.~Agarwal, M.~O'Neil, and M.~Rachh, ``Fmm-accelerated solvers for the
  laplace-beltrami problem on complex surfaces in three dimensions,'' 2021.

\bibitem{DARBAS2006834}
M.~Darbas, ``Generalized combined field integral equations for the iterative
  solution of the three-dimensional maxwell equations,'' \emph{Applied
  Mathematics Letters}, vol.~19, no.~8, pp. 834--839, 2006.

\bibitem{Ciraksubd}
F.~Cirak, M.~Ortiz, and P.~Schr\"{o}der, ``Subdivision surfaces: a new paradigm
  for thin-shell finite-element analysis,'' \emph{Int. J. Numer. Methods Eng},
  vol.~47, pp. 2039--2072, 2000.

\bibitem{subdivision}
D.~Zorin, P.~Schroder, T.~DeRose, L.~Kobbelt, A.~Levin, and W.~Sweldens,
  \emph{Subdivision for Modeling and Animation}.\hskip 1em plus 0.5em minus
  0.4em\relax SIGGRAPH 2000 Course Notes, 2007.

\bibitem{JieDaultShankerChapterSubd}
J.~Li, D.~Dault, and B.~Shanker, \emph{New trends in computational
  electromagnetics}.\hskip 1em plus 0.5em minus 0.4em\relax Institute of
  Engineering and Technology, 2019, ch. New trends in geometric modeling and
  discretization of integral equations, pp. 315--372.

\bibitem{IGA_JCP}
Q.~Pan, T.~Rabczuk, G.~Xu, and C.~Chen, ``Isogeometric analysis for surface
  pdes with extended loop subdivision,'' \emph{Journal of computational
  physics}, vol. 398, p. 108892, 2019.

\bibitem{Rumpf_IGA}
B.~Juettler, A.~Mantzaflaris, R.~Perl, and R.~M, ``On isogeometric subdivision
  methods for pdes on surfaces,'' \emph{Computer methods in applied mechanics
  and engineering}, vol. 302, pp. 131--146, 2016.

\bibitem{Italian_LBO_book}
G.~Patane, \emph{An introduction to Laplacian spectral distances and kernels:
  theory, computation, and applications}.\hskip 1em plus 0.5em minus
  0.4em\relax Morgan \& Claypool, 2017.

\bibitem{Italian_LBO_paper}
------, ``Laplacian spectral basis functions,'' \emph{Computer aided geometric
  design}, vol.~65, pp. 31--47, 2018.

\bibitem{azencot2013operator}
O.~Azencot, M.~Ben-Chen, F.~Chazal, and M.~Ovsjanikov, ``An operator approach
  to tangent vector field processing,'' in \emph{Computer Graphics Forum},
  vol.~32, no.~5.\hskip 1em plus 0.5em minus 0.4em\relax Wiley Online Library,
  2013, pp. 73--82.

\bibitem{spectral_tangential}
\BIBentryALTinterwordspacing
C.~Brandt, L.~Scandolo, E.~Eisemann, and K.~Hildebrandt, ``Spectral processing
  of tangential vector fields,'' vol.~36, no.~6, 2017. [Online]. Available:
  \url{https://doi.org/10.1111/cgf.12942}
\BIBentrySTDinterwordspacing

\bibitem{poelke2017hodge}
K.~Poelke, ``Hodge-type decompositions for piecewise constant vector fields on
  simplicial surfaces and solids with boundary,'' Ph.D. dissertation, 2017.

\bibitem{custers2020subdivision}
B.~Custers and A.~Vaxman, ``Subdivision directional fields,'' \emph{ACM
  Transactions on Graphics (TOG)}, vol.~39, no.~2, pp. 1--20, 2020.

\bibitem{Local_CC_CFIER}
R.~J. Adams, ``Physical and analytical properties of a stabilized electric
  field integral equation,'' \emph{IEEE transactions on antennas and
  propagation}, vol.~52, no.~2, pp. 362--372, 2004.

\bibitem{vikram2009novel}
V.~Melapudi, H.~Huang, B.~Shanker, and T.~Van, ``A novel wideband {FMM} for
  fast integral equation solution of multiscale problems in electromagnetics,''
  \emph{IEEE Trans. Antennas Propag.}, vol.~57, no.~7, pp. 2094--2104, 2009.

\bibitem{dault2016mixed}
D.~Dault and B.~Shanker, ``A mixed potential {MLFMA} for higher order moment
  methods with application to the generalized method of moments,'' \emph{IEEE
  Trans. Antennas Propag.}, vol.~64, no.~2, pp. 650--662, 2016.

\bibitem{song2009modeling}
X.~Song and B.~J{\"u}ttler, ``Modeling and 3d object reconstruction by
  implicitly defined surfaces with sharp features,'' \emph{Computers \&
  Graphics}, vol.~33, no.~3, pp. 321--330, 2009.

\end{thebibliography}

\end{document}